\documentclass[11pt,a4paper]{article}
\usepackage{amsmath, amssymb, amscd}
\usepackage[latin1]{inputenc}
\usepackage[french]{babel}
\usepackage{fancyhdr}
\title{Pincement spectral en courbure de Ricci positive}
\author{J\'er\^ome Bertrand\thanks{Soutenu par la requ\^{e}te 20-101469
du FNRS}\\
Institut de math\'ematiques, \\ Universit\'e de Neuch\^atel, Suisse. \\
E-mail: jerome.bertrand@unine.ch}
\date{\ }
\newtheorem{lem}{Lemme}[section]
\newtheorem{cor}[lem]{Corollaire}
\newtheorem{theo}[lem]{Th{\'e}or{\`e}me}
\newtheorem{prop}[lem]{Proposition}
\newtheorem{defi}[lem]{D{\'e}finition}
\newenvironment{dem}{\bf Preuve : \rm}{\hfill 
$\blacksquare$\par\medbreak}
\DeclareMathOperator{\hes}{Hess}
\DeclareMathOperator{\vol}{\, \rm{vol} \,}
\DeclareMathOperator{\im}{\, \frac{1}{\vol M } \, \int_{M}}
\DeclareMathOperator{\is}{\, \frac{1}{\vol \mathbb{S}^n} \, 
\int_{\mathbb{S}^n}}
\DeclareMathOperator{\ens}{\, {\mathcal{M}}_n \,}

\DeclareMathOperator{\ep}{ \epsilon}
\newcommand{\ov}[1][]{\overline{#1}}
\newtheorem{rem}[lem]{Remarque}

\begin{document}

\maketitle


\begin{center}{ \bf R\'esum\'e}
\end{center}
\begin{quotation}
Dans cet article, nous d\'emontrons que sur les vari\'et\'es riemanniennes 
de dimension $n$ v\'erifiant $Ric \geq (n-1)g$ et pour $k$ dans
$\{1,\ldots,n+1\}$, la $k$\ieme valeur
propre du laplacien est proche de $n$ si et seulement si la vari\'et\'e contient
une partie
Gromov-Hausdorff proche de la sph\`ere $\mathbb{S}^{k-1}$. Pour $k=n+1$,
nous obtenons une nouvelle preuve des r\'esultats de Petersen et Colding 
qui montrent que pour de telles vari\'et\'es, la $(n+1)$\ieme valeur propre 
est proche de $n$  si et seulement si la vari\'et\'e est Gromov-Hausdorff proche
de la sph\`ere 
de dimension $n$.
\end{quotation}
\bigskip

\begin{center}
    {\bf Abstract} 
\end{center}
\begin{quotation}
We show that for $n$-dimensional manifolds with $Ric \geq (n-1)g$ and for $k$ in 
$\{1,\ldots,n+1\}$, the $k$-th eigenvalue for the Laplacian is close to 
$n$ if 
and only if the manifold contains a subset which is Gromov-Hausdorff 
close to the sphere $\mathbb{S}^{k-1}$. For $k=n+1$, this gives a new 
proof of results of Colding and Petersen which show that the $(n+1)$-th 
eigenvalue is close to $n$ if and only if the manifold is 
Gromov-Hausdorff close to the $n$-sphere. 
\end{quotation}
\section*{Introduction}

Dans cet article, nous consid\'erons les vari\'et\'es riemanniennnes connexes, compactes $(M,g)$ de 
dimension $n$ dont la courbure de Ricci v\'erifie l'in\'egalit\'e 
$Ric \geq (n-1)g$. On note $\ens$ l'ensemble (des classes d'isométrie) de ces vari\'et\'es. 
Sur $\ens$, la sph\`ere canonique $(\mathbb{S}^n,g_{can})$ r\'ealise 
l'extr\'emum de plusieurs invariants riemanniens.
\begin{theo}[\cite{myers,bishop,lichnerowicz,cheng1975,obata}] Tout élément $(M,g)$ de $\ens$ vérifie
\begin{center} $\begin{array}{l}
  {\rm diam} \,(M) \leq \pi,  \\
  \vol (M) \leq \vol (\mathbb{S}^n), \\ 
  \lambda_1(M) \geq n, 
\end{array}
$ 
\end{center}
o\`u $\lambda_1(M)$ d\'esigne la premi\`ere valeur propre (non nulle) 
du laplacien de $(M,g)$ agissant sur les fonctions. De plus, dans chaque inégalité, l'égalité n'a lieu que si la variété $(M,g)$ est isométrique
 à la sphère canonique.
 \end{theo}
\smallskip
L'objet de cet article est de caract\'eriser les vari\'et\'es 
appartenant \`a $\ens$ dont le d\'ebut du spectre est presque minimal 
(c'est à dire proche de $n$). De nombreux auteurs se sont intéressés à cette question de la presque égalité (aussi appelée pincement) pour des invariants riemanniens comme le volume ou le diamètre. Nous rappelons ci-dessous quelques-uns de ces résultats, qui ont plus particulièrement motivé ce travail. Le premier d'entre eux est que la premi\`ere valeur propre (non nulle) du laplacien d'un élément de $\ens$ est proche de $n$ si et seulement si le 
diam\`etre de cette vari\'et\'e est proche de $\pi$ (la condition n\'ecessaire est due \`a S.Y. Cheng 
\cite{cheng1975}, la condition suffisante \`a C. Croke 
\cite{croke1982}). Plus précisément, l'équivalence est la suivante :
\begin{theo}[\cite{cheng1975,croke1982}]\label{rajout2} Pour tout réel positif $\ep$, il existe un réel positif $\eta$ tel que tout élément $(M,g)$ de $\ens$ pour lequel ${\rm diam }\, (M)> \pi -\eta$ (respectivement $\lambda_1(M) <n + \eta$), vérifie $  \lambda_1(M) <n + \ep$ (respectivement
$ {\rm diam }\, (M)> \pi-\ep$).\end{theo}
Dans la suite, nous écrirons simplement {\og est proche de \fg} pour ce type d'équivalence.

 En 1996, T. Colding a démontré des résultats de pincement faisant intervenir la distance de Gromov-Hausdorff (nous renvoyons à \cite{gallot1997} pour la définition de cette distance) :
\begin{theo}[\cite{coldingSh1996, coldingLa1996}]\label{a9t4}
Pour tout élément $(M,g)$ de $\ens$, les propriétés suivantes sont équivalentes :

$
\begin{array}{ll}
 1) & \vol(M) \mbox{ est proche de } \vol (\mathbb{S}^n), \\
 2) & {\rm rad}\, (M) \mbox{ est proche de } \pi, \\
 3) & d_{GH}(M,\mathbb{S}^n) \mbox{ est proche de }0,
 \end{array}
$

o\`u {\rm rad}$ \,(M)$ est le plus petit rayon d'une boule recouvrant $M$ et 
$d_{{GH}}$ d\'esigne la distance de Gromov-Hausdorff.
\end{theo}

 Par la suite, P. Petersen a obtenu une nouvelle condition équivalente faisant intervenir le spectre du laplacien. Rappelons que sur la 
sph\`ere 
canonique $\mathbb{S}^n$, la premi\`ere valeur propre non nulle $n$ est de 
multiplicit\'e $n+1$. 
\begin{theo}[\cite{petersen1999}]\label{a9t3}
Pour tout élément $(M,g)$ de $\ens$, les conditions suivantes sont équivalentes :

$\begin{array}{ll}
 1) &{\rm rad }(M) \mbox{ est proche de } \pi, \\
 2) & \lambda_{n+1}(M)  \mbox{ est proche de } n .
 \end{array}
$
\end{theo}

Pour \'enoncer le r\'esultat principal de cet article, nous avons besoin de la 
d\'efinition suivante.
\begin{defi}Soit $(M,g)$ un élément de $\ens$, $k$ un entier 
positif et $\eta$ un nombre r\'eel positif ou nul. La variété 
$(M,g)$ v\'erifie la propri\'et\'e $P_k(\eta)$ s'il existe $k$ 
couples de points $(x_1,y_1),\dots,(x_k,y_k)$ dans $M^2$ v\'erifiant
pour tout $i$ dans $\{1,\ldots,k\}$, 
$$d(x_i,y_i) > \pi  -\eta,$$
et pour tout $i,j$ distincts dans $\{1,\ldots ,k\}$, 
$$\left|d(x_i,x_j)-\frac{\pi}{2}\right| < \eta.$$
\end{defi}

Dans cet article, nous démontrons le
\begin{theo}\label{a9t2}Soit $k$ dans $\{1,\ldots,n+1\}$. Pour tout élément $(M,g)$ de $\ens$, les propriétés suivantes sont équivalentes :

 1)   $ \lambda_k(M)$ est proche de $n $, 
 
 2) $ (M,g)$ v\'erifie  la propriété $P_k(\eta)$ pour $ \eta$ proche de $0$, 
 
 3)  $(M,g)$ contient une partie $A_k$ telle que 
  $d_{{GH}}(A_{k},\mathbb{S}^{{k-1}})$ est proche  
de $0$.

De plus, si la troisième propriété est satisfaite, alors la partie $A_k$ v\'erifie \'egalement 
une propri\'et\'e de {\og presque convexit\'e \fg}. Nous renvoyons \`a la 
proposition \ref{a9p1} pour plus de d\'etails.
\end{theo}

Lorsque $k=1$, l'énoncé du théorème \ref{a9t2} se ramène à celui du théorème \ref{rajout2}. Lorsque $k=n+1$, notre d\'emonstration fournit en particulier une nouvelle preuve de la propriété
\begin{equation}\label{a9e3}
 \lambda_{n+1}(M)\mbox{ proche de } n,  \mbox{ implique } d_{{GH}}(M,\mathbb{S}^n) \mbox{ proche de } 0,
\end{equation}
 autrement dit, on peut prendre $A_{n+1}=M$ dans l'énoncé ci-dessus. 

 Sous les hypoth\`eses du th\'eor\`eme \ref{a9t3}, P. 
Petersen montre \'egalement que l'application $\Phi = 
\frac{(f_1,\ldots,f_{n+1})}{\sqrt{f_1^2+\cdots +f_{n+1}^2}}$ est une 
approximation de Gromov-Hausdorff de $M$ sur la sph\`ere canonique, o\`u les 
$(f_i)_{1\leq i\leq n+1}$ sont les fonctions propres associ\'ees \`a 
$(\lambda_i(M))_{1\leq i\leq n+1}$ et normalis\'ees par analogie avec 
le cas de la sph\`ere. Une \'etape importante de la d\'emonstration de P. Petersen de 
l'implication (\ref{a9e3}) est de prouver que l'application $\Phi$ 
est surjective. Pour cela P. Petersen montre que le degr\'e de $\Phi$ est 
non nul\footnote{P. Petersen a reconnu avoir commis une erreur dans sa démonstration de la surjectivité \cite{petersen04}, la preuve du théorème \ref{a9t2} donne en particulier une preuve complète de la propriété (\ref{a9e3})}. Dans notre cas, il n'y a pas de raison pour que la partie 
$A_k$ soit une vari\'et\'e de dimension $k-1$ et donc on ne peut pas 
appliquer un argument de degr\'e. Nous utilisons \`a la place un 
lemme 
de Toponogov $L^2$ initialement introduit par T. Colding dans 
\cite{coldingSh1996} mais contrairement \`a T. Colding nous ne fixons 
pas des 
conditions au bord  pour l'\'equation diff\'erentielle sous-jacente mais des 
conditions initiales de Cauchy. Ceci nous permet d'obtenir une 
nouvelle 
preuve du r\'esultat de C. Croke et m\^eme d'obtenir un lien pr\'ecis 
entre la fonction propre et les points \`a distance presque $\pi$ 
(voir la proposition \ref{a1l8}). Le contr\^ole de la condition 
initiale sur la d\'eriv\'ee est une cons\'equence d'une estimation du 
gradient d'une fonction propre due \`a P. Li et S.T. Yau.

Cet article est organis\'e de la mani\`ere suivante. Dans la 
deuxi\`eme section, nous donnons les estimations sur les fonctions 
propres 
qui nous seront n\'ecessaires pour d\'emontrer le th\'eor\`eme 
\ref{a9t2}, nous pr\'esentons \'egalement le lemme de Toponogov $L^2$ 
 (lemme \ref{a1l1}) et 
l'utilisation que nous allons en faire. Dans la troisi\`eme section, 
nous d\'emontrons qu'un élément de $\ens$ v\'erifiant la 
propri\'et\'e $P_k(\eta)$ pour $\eta$ petit a n\'ecessairement $k$ 
valeurs propres proches de $n$. La quatri\`eme section est 
consacr\'ee 
\`a la r\'eciproque. Dans la derni\`ere partie nous montrons que la condition
 1) implique la condition 3), la réciproque \'etant imm\'ediate 
 cela termine la preuve du th\'eor\`eme \ref{a9t2}.   
\section{R\'esultats pr\'eliminaires}

Soit $ p \geq 1$ un nombre r\'eel et $h$ appartenant à $L^p (M)$. On note
$$ ||h||_{L^p}= \left(\frac{1}{\vol M} \int_M h^p \,
dx\right)^{\frac{1}{p}}. $$
On utilisera la définition usuelle pour la norme $L^{\infty}$.

{\bf
On notera $\tau (\ep),r(\ep), \eta(\ep)$, etc $\dots$ de mani{\`e}re 
g{\'e}n{\'e}rique, toute \\quantit{\'e} positive ne
d{\'e}pendant que de $\ep$ et de la dimension $n$ de la 
vari{\'e}t{\'e}, dont la limite quand $\ep$ tend vers $0$ est $0$. 

}

Enfin, toute fonction propre $f$ de valeur propre proche de $n$ sera  
normalis\'ee par analogie avec le cas de la sph\`ere, par 
\begin{equation}\label{a9e2}
 \im f^2 = \frac{1}{n+1}.
\end{equation}

\medskip

Les fonctions propres de $(\mathbb{S}^n,g_{can})$ associ\'ees \`a la valeur 
propre $n$ sont les fonctions $\cos d_x$, avec $x$ appartenant à $ \mathbb{S}^n$ et $d_x$ la fonction distance au point $x$. En 
particulier, elles v\'erifient 
\begin{equation}\label{a9e1}
\cos^2d_x + |\nabla \cos d_x|^2=1,
\end{equation}
$$ \frac{1}{\vol (\mathbb{S}^n)}\int_{\mathbb{S}^n} \cos^2 d_x = 
\frac{1}{n+1}.$$
Sur une vari\'et\'e de $\ens$ admettant des valeurs propres proches 
de $n$, on a l'estimation suivante due \`a P. Li :
\begin{prop}[\cite{li1980}]\label{a2l5}\label{a1l3}Soit $(M^n,g)$ 
une vari{\'e}t{\'e} riemannienne compacte de dimension $n$ dont la 
courbure de 
Ricci v\'erifie  $Ric \geq (n-1)g $. Soit $\ov[f]$ une combinaison lin\'eaire de fonctions propres du 
laplacien sur $(M,g)$ 
$$ \ov[f]= \sum_{i=1}^k a_if_i, $$
avec $k$ un entier non nul et $\Delta f_i = \lambda_i f_i$ pour tout 
$i$ dans  $\{1,\ldots,k\}$.\\
Supposons que pour tout $i$ dans $\{1,\ldots,k\}$, 
$$\lambda_i \leq n+ \ep $$
avec $\ep >0$, alors 
$$|| \ov[f]^2 +|d\ov[f]|^2 ||_{L^\infty} \leq (1+ 
\tau(\ep))(n+\ep+1)||\ov[f]||^2_{L^2},$$
o\`u $\tau (\ep)$ est une fonction croissante ne d\'ependant que de 
$n$ et telle que $\lim_{\ep \rightarrow 0} \tau(\ep)= 0$.
\end{prop}
\begin{rem} En particulier si $f$ est une fonction propre de valeur 
propre non nulle $\lambda \leq n+\ep$ et si $\im f^2 =\frac{1}{n+1}$, on 
obtient pour tout $x$ dans $M$,
 $$f^2(x) + |df|^2(x) \leq 1 + \tau(\ep),$$
donc comme $ \im f^2 + |df|^2 = \frac{\lambda +1}{n+1}$, on en 
d\'eduit que (\ref{a9e1}) est {\og stable \fg} pour la norme $L^1$ :
$$ \im |f^2+|df|^2-1| \leq \tau(\ep).$$

D'autre part, si l'on suppose $\ep <1$ alors il existe une constante 
$C(n)$ ne d\'ependant que de la dimension $n$ de $M$, telle que 
\begin{equation}\label{a9e4}
 ||f^2+ |\nabla f|^2||_{L^{\infty}} \leq C(n),
\end{equation}
\noindent nous utiliserons implicitement cette propri\'et\'e par la 
suite.
\end{rem}

La preuve de ce type d'in\'egalit\'e est essentiellement classique, 
elle repose sur une in\'egalit\'e de Sobolev et le proc\'ed\'e 
d'it\'eration 
de Moser. Nous renvoyons \`a \cite{gallot83} et \`a \cite{aubry00}, 
lemme 1.4 pour une d\'emonstration dans ce cas particulier. Une cons\'equence de cette proposition est qu'un élément de $\ens$ admet au plus  $n+1$ valeurs propres proches 
de $n$.
\begin{cor}[\cite{gallot83}]\label{a5c4}Il existe une constante 
$C(n)>0$ telle que pour tout élément $(M,g)$ de  $\ens$,
$$\lambda_{n+2}(M) \geq n+C(n).$$
\end{cor}

\bigskip
Une autre propri\'et\'e caract\'eristique des fonctions $\cos d_x$ 
sur la sph\`ere $(\mathbb{S}^n,g_{can})$ est qu'elles sont solutions 
de 
l'\'equation 
$$ \hes f + f g_{can} = 0.$$
Un r\'esultat de M. Obata \cite{obata} montre que ce sont les 
seules solutions parmi les fonctions r\'eguli\`eres d\'efinies sur 
un élément $(M,g)$ de $\ens$. Cependant, \`a l'aide de la formule de Bochner, on montre (voir par 
exemple \cite{coldingSh1996}, page 178) la 

\begin{prop}\label{a1l4}Il existe une constante $C(n)$ telle que 
tout élément $(M,g)$ de $\ens$ pour lequel $ n \leq \lambda_1(M) \leq n+ \ep,$
v\'erifie l'in\'egalit\'e 

$$||\hes f + fg||_{L^2} \leq C(n)\ep^{\frac{1}{2}} 
||f||_{L^2},$$
o\`u $f$ est une fonction propre associ\'ee \`a $\lambda_1(M)$.
\end{prop}
\begin{rem}Il existe des éléments de $\ens$ non hom\'eomorphes \`a la 
sph\`ere dont la premi\`ere valeur propre est 
arbitrairement proche de $n$ (voir par exemple \cite{anderson1990}). 
Un r\'esultat de S. Gallot (\cite{gallot80}, lemme 3.1) implique 
alors 
qu'on ne peut esp\'erer obtenir une estimation similaire de $\hes 
f+fg$ en norme $L^{\infty}$.
\end{rem} 
Le lemme suivant permet de d\'eduire des 
informations g\'eom\'etriques de cette in\'egalit\'e sur le hessien.

\begin{lem}[\cite{cheeger95}]\label{a1l1}Soit $ (M,g)$ un élément de $\ens$.
 Il existe des constantes ne d{\'e}pendant que de n not\'ees $
  C(n) $ et $ \tilde{C}(n) $ telles que, 
pour $ x_1 $ et $x_2 $ appartenant \`a $ M$, $
  r_1,r_2$ des réels positifs (on note  $ B_i = B(x_i,r_i) $) et pour toute 
fonction continue $h$ sur $M$, on a
\begin{multline} 
 \frac{1}{\vol(B_1 \times B_2)} \int_{ B_1 \times B_2}
\left( \oint_{{\gamma}_{xy}} h^2\right)dxdy
 \leq  \\
 C(n) \left( \frac{1}{\vol B_1} + \frac{1}{\vol 
B_2} \right)
\int_M  h^2(x)dx.
\end{multline}
On obtient en particulier pour $r_1=r_2=r$ \\
$$ \frac{1}{\vol(B_1 \times B_2)} \int_{ B_1 \times B_2}
 \left(\oint_{{\gamma}_{xy}} h^2\right)dxdy
 \leq \frac{\tilde{C}(n)}{V(r)} \frac{1}{ \vol (M)}\int_M
 h^2(x)dx. $$

\end{lem}
\begin{rem}La notation $V(r)$ d\'esigne le volume d'une boule 
g\'eod\'esique de $(\mathbb{S}^n,g_{can})$ de rayon $r$. La notation $ B_1 
\times B_2$ 
d{\'e}signe en r{\'e}alit{\'e} le
sous-ensemble de mesure pleine de ce produit, constitu{\'e} par les
couples admettant une unique g{\'e}od{\'e}sique minimisante les 
reliant (not\'ee $ {\gamma}_{xy}$). La seconde in\'egalit\'e se 
d\'eduit de 
la premi\`ere en utilisant le théorème de Bishop-Gromov (nous renvoyons à \cite{gallot1997} pour un énoncé).
\end{rem}
En appliquant ce lemme \`a la fonction $|\hes f + fg|$ o\`u $f$ est 
la fonction propre associ\'ee \`a $\lambda_1(M)$ normalis\'ee par 
(\ref{a9e2}), dans le cas o\`u $r_1=r_2=r$ et en remarquant que pour 
une g\'eod\'esique param\'etr\'ee
par longueur d'arc 
 $$ |(f\circ \gamma_{xy})''(t) +(f\circ \gamma_{xy})(t)|^2 \leq |\hes f 
+fg|^2, $$
on  d\'eduit de  la proposition \ref{a1l4} 
\begin{equation}\label{a4e3}
  \mbox{$\frac{1}{\vol(B_1 \times B_2)}$} \int_{ B_1 \times B_2}
 \int_0^{d(x,y)} |(f\circ \gamma_{xy})''(t) +(f\circ 
\gamma_{xy})(t)|^2dtdxdy
 \leq 
C(n) \frac{\ep}{V(r)}. 
\end{equation}
Par cons\'equent, pour des rayons $r(\ep)$ convenables ({\it i.e} 
tels que $\frac{\ep}{V(r(\ep))}$ soit petit), l'in\'egalit\'e de 
Byenaim\'e-Tchebitchev implique l'existence de points $x,y$ pour 
lesquels $ f\circ \gamma_{xy}$ v\'erifie presque la m\^eme \'equation 
diff\'erentielle que dans le cas de la sph\`ere canonique. On peut ensuite par des m\'ethodes classiques comparer $ f\circ 
\gamma_{xy}$ \`a une solution correspondante sur la sph\`ere en 
fixant des conditions au bord 
(comme l'a fait T. Colding dans \cite{coldingSh1996}) \`a l'aide du 
lemme suivant.

 \begin{lem}\label{a1l9.5}Soit $v(t)$ et $Z(t)$ deux fonctions 
d{\'e}finies
  sur $ [0,l]$ avec $l<\pi$.  On suppose que $ \int_0^l Z^2(t)dt < 
\epsilon^2 $ et que
 $v$ est solution de $ v"+v=Z $ avec $ |v(0)-a|< \eta $ et   $
 |v(l)-b| < \eta $. Il existe une constante positive $C$ telle que pour tout $t$ dans $[0,l]  $,
 $$ |v(t)-\tilde{u}_{a,b}(t)| < \frac{C}{\sin (l)}(\ep +\eta)  $$
 et
 $$|v'(t)-\tilde{u}_{a,b}'(t)| <\frac{C}{\sin (l)}(\ep +\eta),$$
o\`u $\tilde{u}_{a,b}  $ est la solution de $ {u}" + u=0$ sur $[0,l] 
$ v{\'e}rifiant les conditions initiales $ u(0)=a $ et $ u(l)=b $.
\end{lem}
On peut \'egalement fixer des conditions de Cauchy.
\begin{lem}\label{a1l9}Soit $v(t)$ et $Z(t)$ deux fonctions 
d{\'e}finies
  sur $ [0,l]$ avec $l \leq \pi$. On suppose que $ \int_0^l Z^2(t)dt 
< \epsilon^2 $ et que
 $v$ est solution de $ v"+v=Z $ avec $ |v(0)-a|< \eta $ et  $ 
|v'(0)-b|<\eta $
 . Il existe une constante positive $C$ telle que pour tout $t$ dans  $[0,l], $
 $$ |v(t)-u_{a,b}(t)| < C(\epsilon +\eta)$$
 et 
$$|v'(t)-u_{a,b}'(t)| < C(\epsilon+\eta) ,$$
o\`u $u_{a,b}$ est la solution de $ {u}" + u=0$ sur $[0,l] $ 
v{\'e}rifiant les conditions initiales $ u(0)=a $ et
$ u'(0)=b$.
\end{lem} 

Pour contr\^oler les conditions initiales de l'\'equation 
diff\'erentielle dans le lemme \ref{a1l9}, nous aurons 
besoin d'une estimation due \`a P. Li 
et S.T. Yau (\cite{li1980}, voir \'egalement \cite{schoen94}, page 
108) qui prouve que la norme du gradient d'une fonction propre sur un élément de $\ens$, reste petit au voisinage des points 
r\'ealisant les extr\'ema de la fonction propre. Ce r\'esultat ne 
peut d\'ecouler directement d'une estimation sur le 
hessien de la fonction propre car on constate, en consid\'erant des 
sph\`eres rondes de rayon arbitrairement petit, que la norme 
$L^{\infty}$ 
du hessien d'une fonction propre de norme $1$, tend vers l'infini. 
\begin{prop}[\cite{li1980}]\label{a1p1}
Soit $ (M,g) $ un élément de $\ens$ et $f $ une fonction propre de valeur propre non nulle $ 
{\lambda}  $. Sous ces hypothèses, on a pour tout $x$ dans $M$, l'estimation 
$$  | \nabla f |^2(x) \leq \frac{2\lambda 
\sup_M f}{\sup_M f- \inf_M f}
 (\sup_M f-f(x))(f(x)- \inf_M f) .$$ 
\end{prop}

\section{Vari\'et\'es v\'erifiant la propri\'et\'e $P_k(\eta)$}

L'objet de cette partie est de d\'emontrer le 
\begin{theo}\label{a5t6}Soit $k$ dans $\{2,\ldots,n+1\}$. Il existe une 
fonction $\tau(\eta)$ telle que, pour tout élément  
$(M,g)$ de $\ens $ v\'erifiant la propri\'et\'e $P_k(\eta)$, on a 
l'estimation 

$$ \lambda_k(M) \leq n+ \tau (\eta).$$
\end{theo}
 
La d\'emonstration du th\'eor\`eme \ref{a5t6} repose en partie sur 
l'utilisation de fonctions $\cos d_p$ pour $p$ appartenant \`a $M$ et 
admettant un {\og presque antipode \fg} (c'est à dire $p$ est tel que $\sup_{x \in M} d(p,x)$
est proche de $\pi$). Nous \'etudions de telles fonctions dans le prochain 
paragraphe.
\subsection{Propri{\'e}t{\'e}s des fonctions $\cos d_p$}
Sur la sph\`ere canonique, toute fonction (propre) $\cos d_p$ est une 
combinaison lin\'eaire d'une base de fonctions propres associ\'ees 
\`a la valeur propre $n$. Pr\'ecis\'ement, si $(x_i)_{1\leq i\leq 
n+1}$ est une base orthonorm\'ee de l'espace euclidien 
$\mathbb{R}^{n+1}$ 
alors pour tout élément $p$ de $\mathbb{S}^n$ 

\begin{equation}\label{a5e16}
 \cos d_p = \sum_{i=1}^{n+1} \cos 
d(p,x_i) \cos d_{x_i}.
\end{equation}
En particulier si $p$ appartient \`a $\mathbb{S}^{k-1}$ (en 
identifiant $\mathbb{S}^{k-1}$ \`a la partie de 
$ \mathbb{S}^n$ dont  les $n-k+1$ derni\`eres coordonn\'ees sont 
nulles), seuls les $p$ premiers termes de la somme ci-dessus sont non 
tous 
nuls. Nous allons montrer que la propri\'et\'e (\ref{a5e16}) sur les 
fonctions $\cos d_p$ est {\og stable \fg} pour les points $p$ admettant un 
presque antipode. Ce r\'esultat améliore un lemme d\'emontr\'e par P. Petersen (\cite{petersen1999}, 
lemme 4.3).
\begin{lem}\label{a1l6} 
 Il existe une fonction $\tau(t)$ tendant 
vers $0$ avec $t$, telle que pour tout nombres réels positifs $\ep,\eta$ vérifiant $\eta < \ep$ et pour tout 
élément $(M,g)$ de $\ens$ contenant deux points $p,q$ 
v\'erifiant $d(p,q) > \pi -\eta$, alors on a 

\begin{equation}\label{a1e0.5}
  || \cos d_p - \sum_{i=1}^{k} a_i(p)f_i ||_{L^\infty} \leq 
    \tau (\frac{\eta}{\ep}),
\end{equation}
o{\`u} $ k=max\{i; {\lambda}_i(M) \leq n+  \epsilon \} $ est un 
entier non nul et o\`u 
les $a_i(p)$ sont les coefficients de Fourier de la fonction $\cos 
d_p$ par rapport \`a une base orthogonale $(f_i)_{i \geq 0}$ de 
fonctions propres normalis\'ees par (\ref{a9e2}), c'est \`a dire 
$ a_i(p)= \frac{n+1}{\vol M} \int_M \cos d_p f_i$. De plus, les coefficients $a_i(p)$ v{\'e}rifient pour $\ep $ assez 
petit
$$ |\sum_{i=1}^k {a_i}^2(p) -1| \leq C(n)\frac{\eta}{\ep}.$$  
\end{lem}

\begin{dem}
Soit $p$ comme dans l'\'enonc\'e. Dans la suite, on note $(a_i)_{i 
\in \mathbb{N}}$ les coefficients de Fourier de $\cos d_p$. 
Une cons\'equence de la formule de la coaire et du th\'eor\`eme de 
Bishop-Gromov est le  
\begin{lem}\label{a1l5}  Il existe une constante $C(n)$
  telle que pour tout élément $ (M,g)$ de $\ens, 
$ admettant deux points $p$ et $q$ v{\'e}rifiant $ d(p,q)> \pi 
-\eta$ et pour toute fonction $ u~: [0,\pi] \rightarrow \mathbb{R} $ 
de
  classe $C^1$, on a 
$$ \left| \im u \circ d_p \, dv - \is  u \circ 
d_{\mathbb{S}^n}(\ov[p],.) \, dx \right|
\leq  C(n) \eta \int_0^{\pi} |u'(r)| dr. $$
\end{lem}
\begin{rem} On obtient en particulier que $\frac{1}{\vol (M)} \int_M
\cos^2 d_p$ est proche de $\frac{1}{n+1}$ si $p$ admet un presque antipode.
\end{rem}
 Nous renvoyons \`a \cite{aubry00} pour une d\'emonstration. En appliquant 
 le lemme \ref{a1l5} {\`a} $u= \cos ^2$, $ u=\sin^2$ et
$u=\cos$, on obtient \\
$$ \left| || \nabla \cos d_p ||^2_{L^2} - n|| \cos d_p
  ||^2_{L^2} \right| \leq C(n) \eta, 
$$
$$ |a_0|= \left| \frac{\sqrt{1+n}}{\vol (M)} \int_M \cos d_p  \right| 
\leq C(n) \eta. $$
Par cons\'equent, les coefficients de Fourier de la fonction $\cos d_p$ 
v\'erifient 

$$ \sum_{i=1}^{+ \infty} \lambda_i \,{a_i}^2  ||f_i||^2_{L^2} \leq n
\sum_{i=1}^{+ \infty}{a_i}^2 ||f_i||^2_{L^2} + C'(n) \eta. $$
C'est {\`a} dire, d'apr\`es la normalisation (\ref{a9e2}) des 
fonctions propres
$$\sum_{i=1}^{+ \infty} \lambda_i \,{a_i}^2 \leq n
\sum_{i=1}^{+ \infty}{a_i}^2+ C''(n) \eta. $$

\noindent Comme $\lambda_1(M) \geq n$, on peut n\'egliger les $k$ 
premiers termes, on obtient
$$ \sum_{i=k+1}^{+ \infty} (\lambda_i-n){a_i}^2  \leq C''(n)\eta. $$
On en d{\'e}duit par  d{\'e}finition de $k$,
\begin{equation}\label{a8e3}
  || \cos d_p - \sum_{i=1}^{k} a_if_i ||^2_{L^2} \leq  
C''(n)\frac{\eta}{\ep}
\end{equation}
et comme $\im \cos^2 d_p$ est proche de $\frac{1}{n+1}$, la formule de Parseval
donne 

$$|\sum_{i=1}^k {a_i}^2 -1| \leq \ov[C](n)\frac{\eta}{\ep}.$$
  
L'inégalité précédente montre en particulier que les coefficients 
$(a_i)_{1 \leq i \leq k}$ sont bornés, on déduit alors de l'hypothèse sur $\lambda_k(M)$
et de la proposition \ref{a1l3}, l'existence d'une constante $D(n)$ telle que 
la fonction  $ h=\cos d_p - \sum_{i=1}^{k} a_if_i$ v{\'e}rifie 
$$ ||h||_{L^\infty} \leq D(n)$$
 et 
$$ ||dh||_{L^\infty}\leq D(n).$$
Or par l'in\'egalit\'e de Cauchy-Schwartz, on d\'eduit de 
(\ref{a8e3}) 
$$ ||h||_{L^1} \leq \left(C''(n) 
\frac{\eta}{\ep}\right)^{\frac{1}{2}}. $$
 Soit $x_0$ tel que $ |h(x_0)|= ||h||_{L^\infty} $ et $r$ un réel positif, le théorème des accroissements finis donne 
 $$ \left(C''(n) \frac{\eta}{\ep}\right)^{\frac{1}{2}} \vol (M)\geq 
\int_{B(x_0,r)} |h(x)|dx \geq
   (||h||_{L^\infty} -D(n)r) \vol (B(x_0,r)).$$
En appliquant le th\'eor\`eme de Bishop-Gromov, on en d\'eduit 
$$ \left(C''(n) \frac{\eta}{\ep}\right)^{\frac{1}{2}} \geq 
\frac{V(r)}{V(\pi)} (||h||_{L^\infty}
-D(n)r). $$
D'o{\`u} le r{\'e}sultat en choisissant $r= r(\frac{\eta}{\ep})$
convenable. Ce qui ach\`eve la preuve du lemme \ref{a1l6}.
\end{dem}

\subsection{D\'emonstration du th\'eor\`eme \ref{a5t6}} 
  
Soit $(M,g)$ un élément de $\ens$ v\'erifiant la propri\'et\'e $P_k(\eta)$, en 
particulier $(M,g)$ vérifie 
$$ { \rm diam}(M) \geq \pi - \eta .$$
Sous ces hypoth\`eses, le r\'esultat de S.Y. Cheng cit\'e dans 
l'introduction montre que $\lambda_1(M)$ est proche de $n$. On d\'eduit du 
lemme 
\ref{a1l5} une estimation explicite.
\begin{lem}\label{a5l2}Il existe une constante $C(n)$ telle que  
tout élément $(M,g)$ de $\ens $ pour lequel  ${\rm diam} (M) > \pi -\eta$, v\'erifie 
l'in\'egalit\'e 
$$ \lambda_1(M) \leq n+ C(n) \eta.$$
\end{lem}
\medskip

Notons
$$\ep=\sqrt{\eta}$$
et 
$$ k_{\ep} = \max \, \{k \in \mathbb{N}; \,\lambda_k(M)\leq n+\ep\}.$$
Par le lemme \ref{a5l2}, pour $\eta$ assez petit, $ k_{\ep}$ est supérieure ou égale à $1$. Notons $(f_i)_{1\leq i \leq k_{\ep}}$ une famille orthogonale de 
fonctions propres asoci\'ees \`a $\lambda_i(M)$ et normalis\'ees par 
$\im 
f_i^2 =\frac{1}{n+1}$. D'autre part, notons \\
$$F=Vect_{L^2(M)}\{f_1,\ldots,f_{k_{\epsilon}}\}$$
 et $P_F$ la projection
orthogonale de $L^2(M)$ sur $F$. Par hypoth\`ese sur $(M,g)$, il existe
 $(x_1,y_1),\dots,(x_k,y_k)$ appartenant \`a $M^2$ tels que pour tout $i,j$ distincts dans $\{1,\ldots,k\}$, 
\begin{equation}\label{a9e9}
 |d(x_i,x_j)-\frac{\pi}{2}|\leq \eta 
\end{equation}
et pour tout $i$ dans $\{1,\ldots,k\}$,
$$ d(x_i,y_i) \geq \pi -\eta.$$

\noindent Par cons\'equent, d'apr{\`e}s le lemme \ref{a1l6} 
appliqu\'e avec $\ep=\sqrt{\eta}$ et $\eta$, il existe une fonction 
$\tau(\eta)$ telle que, pour tout $i$ dans $\{1,\ldots,k\}$, 
\begin{equation}\label{a1e7}
  || \cos d_{x_i} - P_F(\cos d_{x_i})||_{L^\infty} \leq \tau(\eta) 
\end{equation}
o\`u $P_F(\cos d_{x_i})= \sum_{j=1}^{k_{\ep}} a_j(x_i)f_j$ v\'erifie pour tout $i$ 
 \begin{equation}\label{a5e30}
  | \sum_{j=1}^{k_{\epsilon}} a_j^2(x_i) 
-1 | \leq \tau(\eta).
\end{equation}
En particulier, pour $\eta$ assez petit et pour tout $i$, $ P_F(\cos d_{x_i})$ n'est pas identiquement nulle. Par conséquent si $k_{\ep} < k$ alors la famille $(P_F(\cos(d_{x_i}))_{i=1}^k$ est 
li{\'e}e. Notons $(b_i)_{i \in \{1,\ldots,k\}} \mbox{ avec } 
\sum_{i=1}^k b_{i}^2 =1$, des coefficients tels que 
$$\sum_{i=1}^k b_{i}P_F(\cos d_{x_i})=0.$$
Alors (\ref{a1e7}) implique
$$ || \sum_{i=1}^k b_i\cos
d_{x_i}-\sum_{i=1}^k b_iP_F(\cos d_{x_i}) ||_{L^\infty} \leq 
k^{\frac{1}{2}}
\tau(\eta).$$
C'est {\`a} dire,
\begin{equation}\label{a1e8}
 || \sum_{i=1}^k b_i\cos d_{x_i} ||_{L^\infty} \leq k^{\frac{1}{2}}
\tau(\eta). 
\end{equation}
 Comme $ \sum_{i=1}^{k} b^2_i=1$, l'un des coefficients $b_{i_0}$ 
v\'erifie  $ |b_{i_0}| \geq \frac{1}{\sqrt{k}} $. Or l'estimation (\ref{a1e8}) appliqu\'ee au point $ x=x_{i_0}$ et 
l'hypoth\`ese (\ref{a9e9}) implique 
$$ |b_{i_0}| \leq \tau(\eta)$$
ce qui est absurde pour $\eta$ assez petit et donc $ k_{\ep} \geq k$, 
ce qui termine la preuve du th\'eor\`eme \ref{a5t6}.

\section{Valeurs propres proches de $n$} 

L'objet de cette partie est de d\'emontrer la r\'eciproque du 
th\'eor\`eme \ref{a5t6}.
\begin{theo}\label{a9t1}Il existe une fonction $\tau(\ep)$ telle que 
tout élément  $(M,g)$ de $\ens$ pour lequel
$\lambda_k(M) \leq n+ \ep $, v\'erifie la propri\'et\'e 
$P_k(\tau(\ep))$.
\end{theo}
\smallskip

\noindent Ce r\'esultat d\'ecoule d'une {\og r\'eciproque \fg} du lemme \ref{a1l6}.
\begin{prop}\label{a1l8}Il existe des fonctions $\tau(\ep)$ et 
$\psi(\ep)$ telles que, pour tout élément  $(M,g)$ de $\ens$ 
et toute fonction propre $f$ sur $M$ de valeur propre non nulle $
  \lambda \leq n+ \epsilon $, normalis\'ee par $\im f^2 = 
\frac{1}{n+1}$, on a l'estimation 
$$ || cos d_x - f||_{L^{\infty}} \leq \tau (\epsilon), $$
avec  $ x$ dans  $M$ tel que $ f(x)= \sup_{M} f $. De plus, si $y$ dans $M$ v{\'e}rifie $f(y)=\inf_{M}f$ 
alors 
$$ d(x,y) > \pi -\psi(\epsilon).$$ 
\end{prop} 

Nous démontrons cette proposition dans le prochain paragraphe.

\subsection{Fonctions propres associ\'ees \`a une {\og petite \fg} valeur 
propre}\label{propre}
La preuve de la proposition \ref{a1l8} est une cons\'equence du lemme \ref{a1l1}.
La premi\`ere \'etape consiste \`a estimer la borne sup\'erieure 
d'une telle fonction propre. 

\begin{lem}\label{a1l20} Il existe une fonction $\tau(\ep)$ ne 
d\'ependant que de $n$, telle que pour tout élément  
$(M,g)$ de $\ens$ et toute fonction propre $f$ sur $M$ de valeur 
propre non nulle $ \lambda \leq
  n+ \ep $, normalis{\'e}e par $ \im f^2 = \frac{1}{n+1}$, on a 
l'estimation 
$$ |\sup_M f-1| \leq  \tau(\ep).$$
\end{lem}

\begin{dem}
Par la minoration de Lichnérowicz de la premi\`ere valeur propre non 
nulle, $\lambda$ v\'erifie 
 $$n \leq \lambda \leq n+ \ep.$$
Par choix de la normalisation de $f$, on a 
\begin{equation}\label{a5e3}
 1 \leq \im f^2 +|df|^2 \leq 1 + \frac{\ep }{n+1} .
\end{equation}
Donc d'apr\`es la proposition \ref{a2l5}, il existe une fonction 
$\tau(\ep)$ telle que $f$ v\'erifie pour tout $x$ dans $M$ 

$$ f^2(x) +|df|^2(x) \leq 1+ \tau(\ep).$$
Ainsi $ f^2 +|df|^2$ est major\'ee par une quantit\'e environ \'egale 
\`a sa moyenne, par cons\'equent $ f^2 +|df|^2$ est $L^1$ proche de 
sa moyenne : 

 $$ \im |1 -f^2(x) -|df|^2(x)| \leq \tau (\ep).$$ 

Soit $x$ v{\'e}rifiant $ f(x)=\sup_M f$. Un corollaire du 
th\'eor\`eme de Bishop-Gromov (\cite{gallot1997}, remarque 2.8) 
implique 
l'existence de $R(\ep)$, $\tau'(\ep)$ 
ne d\'epen\-dant que de $n$ et de $ \tilde{x}$  v\'erifiant 
$d(x,\tilde{x}) \leq R(\ep)$, tels que 
 $$|f^2(\tilde{x})+|df|^2(\tilde {x})-1| \leq \tau'(\ep).$$\\
D'apr\`es la proposition \ref{a1p1}, 
$$  |df|^2(\tilde {x}) \leq \tau(\ep)$$
d'o\`u 
$$ |f(\tilde{x})-1| \leq \tau (\ep),$$
ce qui permet de conclure pour $ f(x)$ puisque le gradient de $f$ est 
born\'e.

\end{dem}

 {\bf Preuve de la proposition \ref{a1l8}}
\smallskip
De l'hypoth\`ese sur la valeur propre, on d\'eduit (proposition 
\ref{a1l4}) 
$$ 
 || \hes f +fg ||_{L^2} \leq \tau(\epsilon). $$
Fixons $x$ comme dans l'\'enonc\'e et soit $u$ dans $M $ quelconque. En 
appliquant le lemme \ref{a1l1} aux boules $B(x,r(\ep))$ et 
$B(u,r(\ep))$ avec $r(\ep)$ convenable, on obtient l'existence d'une 
fonction $\tau(\ep)$ telle que pour tout $u$ dans $M$, il existe $\tilde{u},\tilde{x}$ dans  $M$ 
tels que :\\
- il existe une unique g\'eod\'esique minimisante $\gamma$ reliant 
$\tilde{x}$ \`a $\tilde{u}$,\\
- $d(u,\tilde{u})\leq r(\ep)$, $d(x,\tilde{x})\leq r(\ep)$ et
\begin{equation}\label{a5e6}
\int_0^{d(\tilde{x},\tilde{u})} |(f\circ \gamma)''(t) +(f\circ 
\gamma)(t)|^2dt \leq \tau(\ep).
\end{equation}

Nous allons maintenant estimer les conditions initiales $f\circ 
\gamma (0)$ et $(f\circ \gamma)'(0)$ afin d'appliquer le lemme  
\ref{a1l9}. Par le lemme \ref{a1l20}, on a 
$$ |\sup_Mf-1| \leq \tau_2(\ep).$$

\noindent Par cons\'equent comme $\tilde{x}$ est proche de $x$ et que le 
gradient de $f$ est born\'e, on en d\'eduit l'existence d'une 
fonction 
$\tau_3(\ep)$ telle que 
\begin{equation}\label{a5e4}
| (f\circ \gamma) (0) -1 | \leq \tau_3(\ep) .
\end{equation}
D'autre part, par la proposition \ref{a1p1}, il existe une fonction 
$\tau_4(\ep)$ telle que 
$ |\nabla f| (\tilde{x}) \leq \tau_4(\ep)$ d'o\`u  
\begin{equation}\label{a5e5}
| (f\circ \gamma)'(0)| \leq \tau_4(\ep).
\end{equation}
Gr\^ace \`a (\ref{a5e6}), (\ref{a5e4}) et (\ref{a5e5}), le lemme 
\ref{a1l9} appliqu\'e \`a $f\circ \gamma $ et $\cos$,
 implique l'existence d'une fonction 
$\tau_5(\ep)$ telle que pour tout $t$ dans $[0,d(\tilde{x},\tilde{u})]$,
\begin{equation}\label{a5e7}
 | (f\circ \gamma) (t) - 
\cos (t) | \leq \tau_5(\ep).
\end{equation}
\begin{equation}\label{a5e8}
 | (f\circ \gamma)' (t) 
+ \sin (t) | \leq \tau_5(\ep).
\end{equation}
 En particulier,
$$ |f(\tilde{u})- \cos (d(\tilde{x},\tilde{u}))| \leq \tau_5(\ep).$$
Donc par construction de $\tilde{x}$ et $\tilde{u}$ et comme le 
gradient de $f$ est born\'e, on en d\'eduit l'existence d'une 
fonction 
$\tau_6(\ep)$ telle que 
$$|| f - \cos d(x,.) ||_{L^\infty} \leq \tau_6 (\ep).$$

Montrons maintenant la deuxi\`eme partie de l'\'enonc\'e. Soit $y$ v{\'e}rifiant 
$ f(y)=\inf_{M} f$ et soit $\tilde{x}$ et 
$\tilde{y}$ comme 
 ci-dessus.\\
D'apr\`es (\ref{a5e8}) 
\begin{equation}\label{a9e5}
 |(f\circ \gamma)'(d(\tilde{x},\tilde{y})) + \sin 
(d(\tilde{x},\tilde{y}))| \leq \tau_5(\ep).
\end{equation}
Comme $\tilde{y}$ est proche de $y$ qui est un point r\'ealisant le 
minimum de $f$, on d\'eduit de la proposition 
\ref{a1p1} apliqu\'ee \`a $\tilde{y}$ et de (\ref{a9e5}) 

$$| \sin (d(\tilde{x},\tilde{y}))| \leq \tau_7(\ep).$$
La borne sur le gradient de $f$ et le lemme \ref{a1l20} excluent l'hypothèse que $d(\tilde{x},\tilde{y})$ soit proche de $0$ donc les points $x$ et $y$ sont 
n\'ecessairement \`a distance presque $\pi$.

\begin{flushright}
 $\blacksquare$
 \end{flushright}

\subsection{ D\'emonstration du th\'eor\`eme \ref{a9t1}}
 
Pour cela, on prouve un r\'esultat un peu plus pr\'ecis :
\begin{prop}\label{a1p2} Soit $k$ dans $\{2,\ldots,n+1\}$. Il existe une 
fonction $\tau(\ep)$ telle que pour tout élément  
$(M,g)$ de $\ens$ v\'erifiant  $ {\lambda}_k(M) \leq n+ \epsilon $, 
il existe $(x_1,y_1),\dots,(x_k,y_k)$ dans $M^2$ tels que, 
pour tout $i$ dans  $\{1,\ldots ,k\}$,
 $$ \vert d(x_i,y_i)- \pi \vert \leq \tau   (\epsilon),$$
 pour tout $i,j$ distincts dans $\{ 1,\ldots,k \}$
 $$ \vert d(x_i,x_j)- \frac{\pi}{2} \vert \leq   \tau (\epsilon),$$
 $$  \vert d(x_i,y_j)- \frac{\pi}{2} \vert \leq   \tau (\epsilon),$$
 $$  \vert d(y_i,y_j)- \frac{\pi}{2} \vert \leq   \tau (\epsilon).$$ 
De plus, ces points v{\'e}rifient  $ f_i(x_i)=\sup_M f_i$ et $ 
f_i(y_i)=\inf_M f_i$ avec $(f_i)_{1 \leq i \leq k}$ une famille 
orthogonale 
de fonctions propres associ\'ees \`a \\ $(\lambda_i(M))_{1 \leq i \leq 
k}$ et  normalis\'ees par (\ref{a9e2}).
\end{prop}

Ces couples de points correspondent dans le cas mod\`ele, aux couples 
form\'es de $k$ vecteurs $(e_i)_{1 \leq i \leq k}$ de la base 
canonique de $\mathbb{R}^{n+1}$ et des $k$ vecteurs oppos\'es 
$(-e_i)_{1 \leq i \leq k}$.

 \medskip

\begin{dem}Soit $(x_i)_{1 \leq i \leq k}$ et $(y_i)_{1 \leq i \leq 
k}$ d\'efinis par 
$$  f_i(x_i)=\sup_M f_i \mbox{ et }  f_i(y_i)=\inf_M f_i,$$
avec les fonctions $(f_i)_{1 \leq i \leq k}$ d\'efinies comme ci-dessus. D'apr{\`e}s  la proposition \ref{a1l8}, pour tout $i$ dans $\{1,\ldots ,k\}$, on a
\begin{equation}\label{a9e8} 
 d(x_i,y_i) > \pi -\psi (\epsilon).
\end{equation}
Nous allons montrer  l'existence d'une fonction $\tau(\ep)$, telle 
que pour tout $i,j$ comme dans l'énoncé, 
\begin{equation}\label{a5e13}
  \vert d(x_i,x_j)- \frac{\pi}{2} \vert \leq   \tau (\epsilon).
\end{equation}
Admettons provisoirement ce r\'esultat, on en d\'eduit les autres 
estimations \`a l'aide d'un lemme sur la fonction {\og excess \fg}, d\^u 
\`a K. 
Grove et P. Petersen.
\begin{lem}[\cite{grove90}]\label{a1l2} Il existe une fonction $
  \tau(\epsilon) $ telle que pour tout élément  
$(M,g)$ de $\ens$ et pour tout $p,q$ dans $M $ v{\'e}rifiant  
$$ d(p,q) \geq \pi - \epsilon, $$
 on a 
$$ e_{p,q}(M)=\sup_{x \in M}(d(p,x)+d(q,x)-d(p,q)) \leq  
\tau(\epsilon). $$
\end{lem}
D'apr\`es ce lemme et (\ref{a9e8}), il existe une fonction 
$\tau_2(\ep)$, telle que pour tout $i,j$ distincts dans $\{ 1,\ldots,k \}$, 
$$|d(x_j,y_j)-d(x_i,x_j) -d(x_i,y_j)| \leq \tau_2(\ep),$$
donc par (\ref{a9e8}) et (\ref{a5e13}), il existe une fonction 
$\tau_3(\ep)$ telle que 
$$  |d(x_i,y_j)-\frac{\pi}{2}| \leq \tau_3(\ep).$$
On d\'eduit de mani\`ere similaire l'estimation sur $d(y_i,y_j)$.

D\'emontrons maintenant l'estimation (\ref{a5e13}). Fixons $i,j$ distincts dans \\$\{ 1,\ldots,k \}$. Notons
 $$h= f_if_j + <df_i,df_j>.$$
 Par hypoth\`ese sur les fonctions $(f_i)_{1\leq i \leq k}$ 
$$ \im h =0.$$
Calculons la diff\'erentielle de $h$ 
$$ dh = f_jdf_i +f_i df_j + \hes f_i (\nabla f_j,.) + \hes f_j(\nabla 
f_i,.),$$
$$ dh= (\hes f_i +f_ig)(\nabla f_j,.)+ (\hes f_j +f_jg)(\nabla 
f_i,.).$$
Donc 
$$ |dh|^2 \leq 2\left(||\hes f_i+f_ig||^2\times ||\nabla f_j||^2 + 
||\hes f_j+f_jg||^2\times ||\nabla f_i||^2\right).$$
Or par la proposition \ref{a1l4}, il existe une fonction 
$\tau_4(\ep)$ telle que pour tout $s$ dans $\{1,\ldots,k\}$,
$$||\hes f_s + f_sg||_{L^2} \leq 
\tau_4 (\ep) .$$  

 \noindent Par (\ref{a9e4}), on en d\'eduit l'existence d'une fonction 
$\tau_5(\ep)$ telle que 
$$ \im |dh|^2  \leq \tau_5(\ep).$$

\noindent En appliquant l'in\'egalit\'e de Poincar\'e \`a la fonction $h$ de 
moyenne nulle, on obtient, puisque $\lambda_1(M) \geq n$, l'existence 
d'une fonction $\tau_6(\ep) $ telle que 
$$ \im h^2 \leq \tau_6(\ep).$$ 
On en d\'eduit alors par un corollaire du th\'eor\`eme de 
Bishop-Gromov (\cite{gallot1997}, remarque 2.8) l'existence de fonctions  $\tau_7(\ep)$ et $R(\ep)$ telles que pour tout $x$ dans $M$, il existe $\tilde{x}$  tel que $ d(x,\tilde{x}) < R(\ep)$ et 
 $$| h(\tilde{x})| \leq   \tau_7(\ep).$$
 Pour $x=x_j$, la proposition \ref{a1p1} implique l'existence d'une constante $C'(n)$ telle que 
$$|\nabla f_j|(\tilde{x}_j) \leq C'(n) R(\ep)$$
puisque $f_j(x_j)= \sup_M f_j$. On en d\'eduit qu'il existe une fonction $\tau_8(\ep)$ telle que 
$$ | h(\tilde{x}_j)-f_i(\tilde{x}_j)f_j(\tilde{x}_j)| \leq 
\tau_8(\ep).$$
Mais, par la proposition \ref{a1l8} 
$$|| \cos d_{x_j} - f_j||_{L^{\infty}} \leq \tau_9(\epsilon),$$
d'o\`u 
$$ |f_j(\tilde{x}_j)-1| \leq C(n)R(\ep)+\tau_9(\epsilon).$$
Par cons\'equent, il existe une fonction $\tau_{10}(\ep)$ telle que 
$$ |f_i(\tilde{x}_j)| \leq \tau_{10}(\ep).$$
On termine la preuve en appliquant de nouveau la proposition 
\ref{a1l8}.

\end{dem}

\section{ Proximit\'e de Gromov-Hausdorff}
Dans cette partie, nous montrons le
\begin{theo}\label{a1t2} Soit $ k$ dans $\{ 2,\ldots,n+1 \} $. Il existe une 
fonction $ \tau
  (\epsilon) $ telle que tout élément  $ (M,g)$ 
de $\ens$ vérifiant $\lambda_k(M) \leq 
n+\ep$, poss\`ede
un sous-ensemble $ A \subset M $ presque convexe tel que  
$d_{GH}(A,S^{k-1})\leq \tau(\epsilon)$.
\end{theo}

Nous renvoyons \`a la proposition \ref{a9p1} pour la 
d\'efinition de la presque convexit\'e. Fixons $k$ dans $\{2,\ldots,n+1\}$. Notons $F=(f_1,\ldots,f_k)$ et 
$\Phi=\frac{F}{\sqrt{f_1^2+\dots+f_k^2}}$ o\`u les fonctions $(f_i)_{1 \leq i \leq 
k}$ sont les fonctions propres 
associ\'ees aux valeurs propres $(\lambda_i(M))_{1 \leq i \leq k}$ et 
normalis\'ees par (\ref{a9e2}). Sur la sphère canonique, les fonctions coordonn\'ees (qui 
forment une base de fonctions propres de valeur propre $n$) fournissent un 
plongement isom\'etrique d'une partie de $(\mathbb{S}^n,g_{can})$ sur 
$\mathbb{S}^{k-1} \subset \mathbb{R}^k$ :
$$\begin{array}{ccc}
\{x \in \mathbb{S}^n; X_1^2+\dots+X_k^2=1\} & \longrightarrow & 
\mathbb{S}^{k-1} \\
 x & \longmapsto & (X_1,\ldots,X_k)(x) 
\end{array}
$$

Nous allons montrer que l'application $\Phi$ restreinte \`a une 
partie $A_k$ convenable est une $\tau(\ep)$ approximation  de 
Gromov-Hausdorff, c'est \`a dire que pour tout $X$ dans  $\mathbb{S}^{k-1}$, il existe $y$ dans
 $A_k$ tel que $d_{\mathbb{S}^{k-1}}(X,\Phi(y)) < \tau(\ep)$
et pour tout  $x,y$ dans $A_k$, $|d(x,y)-d_{\mathbb{S}^{k-1}}(\Phi(x),\Phi(y))| < 
\tau(\ep)$ ($d_{\mathbb{S}^{k-1}}$ désigne la distance induite par la métrique canonique). Nous montrerons que pour une fonction $\eta(\ep)$ bien choisie, 
$A_k=\{x \in M ; |(f_1^2+\dots+f_k^2)(x)-1| < \eta(\ep)\}$ convient.

Par choix de la partie $A_k$ et par uniforme continuit\'e de la 
fonction $\arccos$, il suffit pour d\'emontrer le th\'eor\`eme 
\ref{a1t2} de prouver l'existence d'une fonction $\tau(\ep)$ pour 
laquelle la fonction $F$ vérifie les propriétés suivantes :

\noindent - une propri\'et\'e de {\og $\tau(\ep)$-presque surjectivit\'e \fg} :\\
 pour tout $X$ dans $\mathbb{S}^{k-1}$, il existe $x$ dans $A_k$ tel que
\begin{equation}\label{a5e21} 
||F(x)-X||_{\mathbb{R}^{k}} \leq
  \tau(\epsilon),
\end{equation}
\noindent - une propri\'et\'e de {\og $\tau(\ep)$-proximit\'e m\'etrique\fg} :\\
 pour tout $x,y$ dans $A_k$ 
\begin{equation}\label{a5e22}
| \cos d(x,y) -<F(x),F(y)>_{\mathbb{R}^k}| 
\leq \tau(\ep),
\end{equation}
où $<.,.>_{\mathbb{R}^k}$ désigne le produit scalaire euclidien dans $\mathbb{R}^k$ et $||.||_{\mathbb{R}^{k}}$ la norme associée.

La d\'emonstration de (\ref{a5e21}) fait l'objet du prochain 
paragraphe, nous en d\'eduirons (\ref{a5e22}) dans le paragraphe 
suivant.

\subsection{D\'emonstration de la {\og presque surjectivit\'e\fg}}

Soit $(M,g)$ un élément de $\ens$, $s$ un entier non nul et $\eta$ un nombre r\'eel positif. Soit
 $f_1,\ldots,f_s$ une famille orthogonale de fonctions propres sur 
$M$ associ\'ees \`a $(\lambda_i(M))_{1 \leq i\leq s}$ et 
normalis\'ees 
par (\ref{a9e2}). On note 
$$ A^\eta_s = \{ x \in M ; \; |f_1^2(x)+\dots+f_s^2(x)-1| < \eta\}.$$ 
\begin{prop}\label{a1p3} Soit $ k$ dans $\{1,\ldots,n+1\} $. Il existe des 
fonctions $\eta(\ep)$ et $\tau (\ep)$  telles que pour tout 
élément $(M,g)$ de $\ens$ pour lequel  $\lambda_k(M) \leq
  n+\epsilon $, l'ensemble $A^{\eta(\ep)}_k $ v\'erifie une propri\'et\'e de $\tau(\ep)$-presque surjectivit\'e. 
\end{prop}

\begin{dem}
La preuve repose sur une r\'ecurrence finie. Lorsque $k=1$, la proposition \ref{a1p3} est une conséquence directe de la proposition \ref{a1l8}. Fixons $k$ dans $\{1,\ldots,n+1\} $ et $m$ dans $\{1,\ldots,k-1\}$. Dans la suite de la démonstration, on identifie $\{x \in \mathbb{S}^{k-1}; 
x=(a_1,\ldots,a_m,0,\ldots,0) \}$ avec $\mathbb{S}^{m-1}$.
 La preuve de la proposition est une conséquence du lemme suivant
\begin{lem}\label{a5l1}
Supposons qu'il existe des fonctions $\eta_m(\ep)$ et $\psi_m(\ep)$ 
telles que pour tout $X$ dans $\mathbb{S}^{m-1}$, il existe $x$ dans $A^{\eta_m(\ep)}_m$ tel que
$$  ||(f_1,\ldots,f_m)(x)-X||_{\mathbb{R}^{m}} \leq 
\psi_m(\ep)$$
alors il existe des fonctions  $\eta_{m+1}(\ep)$ et $\psi_{m+1}(\ep)$ 
telles que pour tout $X$ dans $\mathbb{S}^{m}$, il existe $x$ dans $A^{\eta_{m+1}(\ep)}_{m+1}$ tel que
$$ ||(f_1,\ldots,f_{m+1})(x)-X||_{\mathbb{R}^{m+1}} 
\leq \psi_{m+1}(\ep).$$
\end{lem}

Démontrons le lemme \ref{a5l1}. Commen\c{c}ons par remarquer que pour toute fonction $\eta(\ep)$, il existe une fonction 
$\theta(\ep)$ telle que pour tout $s$ dans $\{1,\ldots,k-1\}$, 

\begin{equation}\label{a5e12}
  A^{\eta(\ep)}_s \subset  A^{\theta(\ep)}_{s+1},
\end{equation}
c'est \`a dire que $f_{s+1}(A^{\eta(\ep)}_s)$ est presque réduit à $\{0\}$. Ce r\'esultat est une cons\'equence directe d'un lemme 
d\'emontr\'e par P. Petersen (\cite{petersen1999}, lemme 3.3).
 \begin{lem}[\cite{petersen1999}]\label{a1l11}Il existe une fonction 
$ \tau (\epsilon) $ telle que pour tout élément   
$(M,g)$ dans $\ens $ v\'erifiant  $ \lambda_k(M) \leq n+ \epsilon $, on 
a pour tout $x$ dans $M$   
$$ f_1^2(x)+\dots+f_k^2(x) \leq 1+ \tau 
(\epsilon), $$
o\`u $(f_i)_{1\leq i\leq k}$ est une famille orthogonale de fonctions 
propres de $M$, de valeurs propres $(\lambda_i(M))_{1\leq i\leq k}$ 
et normalis\'ees par $  \im f_i^2 =\frac{1}{n+1}.$
\end{lem}

Soit $ Y=(\cos s_1,\dots,\cos s_{m+1})$ dans $\mathbb{S}^{m} $. 
Il faut distinguer les cas \\$|\sin s_{m+1}|<\mu(\ep)$ et 
$|\sin s_{m+1}| \geq \mu(\ep)$, o\`u $\mu(\ep)$ v\'erifie $  
\lim_{\ep \rightarrow 0} \mu (\ep) =0$ et  sera d\'efini plus loin.

\smallskip

Supposons $|\sin s_{m+1}|<\mu(\ep)$. Dans ce cas, comme $d(x_i,x_{m+1})$ est proche de $\frac{\pi}{2}$ et $d(x_i,y_{m+1})$ est proche de 
$\frac{\pi}{2}$ pour tout $i$ dans  $\{1,\ldots,m\}$, la proposition \ref{a1l8} implique l'existence
 d'une fonction $\tau_2(\ep)$ telle que 
$$ ||(f_1,\ldots,f_{m+1})(\alpha)-Y||_{\mathbb{R}^{m+1}} \leq 
\tau_2(\ep),$$
avec $\alpha=x_{m+1}$ si $s_{m+1}$ est proche de $0$ et $\alpha=y_{m+1}$ sinon.

\smallskip

On suppose maintenant que   $|\sin s_{m+1}| \geq \mu(\ep)$ et que $s_{m+1} 
\leq \frac{\pi}{2}$. Le cas $s_{m+1} > \frac{\pi}{2}$ sera trait\'e plus loin. On d{\'e}finit  
$$X= \left(\frac{\cos s_1}{\sin s_{m+1}},\ldots,\frac{\cos s_m}{\sin s_{m+1}}\right) 
\in \mathbb{S}^{m-1}.$$

\noindent Par hypoth{\`e}se de r{\'e}currence, il existe $x_0 \in 
A^{\eta_m(\ep)}_m  $ telle que 
\begin{equation}\label{a1e6.5}
||(f_1,\ldots,f_m)(x_0)-X||_{\mathbb{R}^m} \leq \psi_m(\epsilon). 
\end{equation}
D'apr\`es (\ref{a5e12}), il existe une fonction $\theta(\ep)$ telle 
que $A^{\eta_m(\ep)}_m \subset A^{\theta(\ep)}_{m+1}$. D'apr\`es la proposition \ref{a1l8}, $f_{m+1}$ est proche en norme $L^{\infty}$ de $\cos 
d_{x_{m+1}}$ donc comme $x_0$ appartient à $A^{\eta_m(\ep)}_m$, l'équation (\ref{a5e12}) 
implique l'existence d'une fonction $\tau_3(\ep)$ telle que 
 $$ |\cos d_{x_{m+1}}(x_0)| \leq \tau_3(\ep).$$
 Par cons{\'e}quent, par l'in\'egalit\'e des accroissements finis, on a 
\begin{equation}\label{a1e6.59}
  |\frac{\pi}{2} - d_{x_{m+1}}(x_0)| \leq \tau_3(\epsilon). 
\end{equation}

En appliquant le lemme \ref{a1l1} aux fonctions $f_i$ (pour 
$i$ dans $\{1,\ldots,m+1\}$) au voisinage des points $x_0$ et $x_{m+1}$, on 
en déduit qu'il existe  $\tilde{x}_0$ avec 
$d(\tilde{x}_0,x_0)<r(\epsilon)$ et $\tilde{x}_{m+1}$  avec $
d(x_{m+1},\tilde{x}_{m+1}) < r(\epsilon)$, tels que si on note 
$\gamma$
l'unique g{\'e}od{\'e}sique minimisante reliant $ \tilde{x}_{m+1}$ 
{\`a} $\tilde{x}_0$ et
$u_i =f_i \circ \gamma $ (pour $i \in \{1,\ldots,m+1\} $) alors 
$$ \int_0^{d(\tilde{x}_0,\tilde{x}_{m+1})} |u_i''(t)+ u_i(t)|^2 dt 
\leq \tau_{4}(\ep).$$
 Pour $i$ dans $\{1,\ldots,m\}$, les conditions aux bords 
sont 
\begin{equation}\label{a8e1}
 u_i(0)=f_i(\tilde{x}_{m+1})
\end{equation}
et
\begin{equation}\label{a8e2} 
u_i(d(\tilde{x}_{m+1},\tilde{x}_0))=f_i(\tilde{x}_0).
\end{equation}

\noindent Or, d'une part 
$$  |f_i(\tilde{x}_{m+1})|\leq |f_i(\tilde{x}_{m+1})-f_i(x_{m+1})| + 
||f_i-\cos d_{x_i}||_{L^{\infty}} + |\cos d_{x_i}(x_{m+1})|, $$
donc
$$ |f_i(\tilde{x}_{m+1})| \leq C(n)r(\epsilon) +
\tau_5(\epsilon) + \tau_6(\epsilon) $$
et d'autre part 
$$ \left|f_i(\tilde{x}_0)- \frac{\cos s_i}{\sin s_{m+1}}\right| \leq 
|f_i(\tilde{x}_0)-f_i(x_0)| + \left|f_i(x_0) - \frac{\cos s_i}{\sin 
s_{m+1}}\right| $$
donc
$$ |f_i(\tilde{x}_0)- \frac{\cos s_i}{\sin s_{m+1}}| \leq
C(n)r(\epsilon) + \psi_m(\ep). $$
Pour $i$ dans $\{1,\ldots,m\}$, notons $\ov[u]_i(t)= \frac{\cos s_i}{\sin 
s_{m+1}} \sin (t)$. En utilisant (\ref{a1e6.59}), on en d\'eduit 
$$|\ov[u]_i(\frac{\pi}{2})-\ov[u]_i(d(\tilde{x}_{m+1},\tilde{x}_0))| 
\leq
\frac{\tau_3(\ep)+2r(\ep)}{\mu(\ep)}.$$
On fixe $\mu(\ep) = \sqrt{\tau_3(\ep)+2r(\ep)}$. D'apr\`es (\ref{a8e1}) et (\ref{a8e2}), il existe une fonction 
$\tau_7(\ep)$ telle que 
$$ |u_i(0) -\ov[u_i](0)| \leq \tau_7(\ep)$$
et 
$$ 
|u_i(d(\tilde{x}_{m+1},\tilde{x}_0))-\ov[u_i](d(\tilde{x}_{m+1},\tilde{x}_0))| 
\leq \tau_7(\ep).$$
D'apr{\`e}s (\ref{a1e6.59}), on peut supposer $\ep$ assez petit pour 
que
$l=d(\tilde{x}_{m+1},\tilde{x}_0)$ v{\'e}rifie l'hypoth{\`e}se du 
lemme
\ref{a1l9.5},  par cons{\'e}quent en appliquant ce lemme  aux 
fonctions
$u_i$ et $\ov[u]_i$, on obtient l'existence d'une fonction 
$\tau_8(\ep)$ telle que pour tout $i$ dans  $\{1,\ldots,m\}$ et pour tout $t$ dans  
$[0,d(\tilde{x}_{m+1},\tilde{x}_0)]$ 
$$ | (f_i \circ \gamma) (t) - 
\frac{\cos
  s_i}{\sin s_{m+1}}\sin (t)| \leq \tau_8(\epsilon). $$
La proposition \ref{a1l8} permet  d'estimer $ f_{m+1}$, on obtient
$$ |f_{m+1}(\gamma(t)) - \cos t | \leq 
\tau_5(\epsilon)+r(\ep).$$
En combinant ces r{\'e}sultats, on en déduit l'existence d'une
fonction $ \tau_9(\epsilon)$ telle que pour tout $t$ dans $[0, d(\tilde{x}_{m+1},\tilde{x}_0)]$,
\begin{multline}
 || (f_1,\ldots,f_{m+1})(\gamma (t)) \\-\left( (   \cos s_1,\ldots,\cos
s_m,0)\frac{\sin t}{\sin s_{m+1}}+(0,\ldots,0,\cos 
t)\right)||_{\mathbb{R}^{m+1}} 
\leq \tau_9(\epsilon).
\end{multline}

\noindent D'o\`u, comme on suppose $s_{m+1} \leq \frac{\pi}{2}$, on 
obtient pour $$T=\min\{d(\tilde{x}_{m+1},\tilde{x}_0),s_{m+1}\}$$  
\begin{multline*}
 || (f_1,\ldots,f_{m+1})(\gamma (T)) - (\cos s_1,\ldots,\cos s_m,\cos 
s_{m+1})||_{\mathbb{R}^{m+1}} \\
\leq  \tau_9(\ep) + \left|\left|\left( (\cos s_1,\ldots,\cos 
s_m,0)\left(\frac{\sin T}{\sin s_{m+1}}-1 \right) + 
\right.\right.\right.\\
\left.\left.\left.(0,\ldots,0,1)(\cos T -\cos 
s_{m+1})\right)\right|\right|_{\mathbb{R}^{m+1}}.
\end{multline*}

\noindent Or comme par (\ref{a1e6.59}) 
 $$ d(\tilde{x}_{m+1},\tilde{x}_0) \geq 
\frac{\pi}{2}-\tau_3(\epsilon)-2r(\epsilon),$$
la proposition est d\'emontr\'ee dans ce cas. Si $s_{m+1} \geq \frac{\pi}{2}$,  il suffit de remplacer $x_{m+1}$ 
par $y_{m+1}$ dans tout ce qui pr{\'e}c{\`e}de.

\end{dem}

\subsection{Propri\'et\'e de {\og proximit\'e m\'etrique\fg}}
Soit $k$ dans $\{2,\ldots,n+1\}$. On conserve la notation $A^{\eta(\ep)}_k$ 
pour la partie introduite dans la proposition \ref{a1p3}. Dans cette partie, nous montrons que $F=(f_1,\ldots,f_k)$ v\'erifie la 
propri\'et\'e de $\tau(\ep)$-proximit\'e m\'etrique, dont nous rappelons la définition :

 Il existe une fonction $ \tau(\ep)$ telle que tout 
élément $ (M,g)$ de $\ens$ tel que $\lambda_k(M) \leq n+\ep$, vérifie
pour tout $x,y$ dans $A^{\eta(\ep)}_k$
\begin{equation}\label{a5e17}
| \cos d(x,y) - 
<F(x),F(y)>_{\mathbb{R}^k}| \leq \tau(\ep),
\end{equation}
ce  qui termine la preuve du th\'eor\`eme \ref{a1t2}. Nous démontrons également que tout élément 
$(M,g)$ de $\ens$ tel que $\lambda_{n+1}(M) \leq n+\ep$, vérifie
\begin{equation}\label{a5e19}
d_{GH}(M,\mathbb{S}^n) \leq \tau (\ep).
\end{equation}
 Enfin nous d\'emontrons que la partie $A^{\eta(\ep)}_k$ est {\og presque convexe \fg} au sens de la proposition ci-dessous. 
\begin{prop}\label{a9p1}Il existe des fonctions $\eta'(\ep)$  et 
$\tau(\ep)$, telles que pour tout élément  $(M,g)$ dans $\ens $
 v\'erifiant $\lambda_k(M)\leq n+\ep$, alors pour tout $x,y$ dans  $A^{\eta(\ep)}_k$, $$d_{A^{\eta'(\ep)}_k}(x,y) \leq d(x,y) + \tau(\ep)$$
o\`u $d_{A^{\eta'(\ep)}_k}$ d\'esigne la distance intrins\`eque de 
l'ouvert $A^{\eta'(\ep)}_k$ et $\eta'(\ep)$ v\'erifie $\eta'(\ep) \geq 
\eta (\ep)$.
\end{prop}
Commen\c{c}ons par montrer comment la démonstration de la propriété (\ref{a5e19}) se ram\`ene \`a la 
d\'emonstration que nous allons donner de la propriété (\ref{a5e17}). Sous l'hypothèse  $\lambda_{n+1}(M) \leq n+\ep$, la proposition 
\ref{a1p3} montre que la variété $M$ 
contient une partie qui est $\tau(\ep)$-presque 
surjective sur $\mathbb{S}^n$. Il suffit donc de prouver que cette 
application 
v\'erifie la propri\'et\'e de $\tau(\ep)$-proximit\'e m\'etrique pour une 
fonction $\tau(\ep)$ convenable. Or d'apr\`es un r\'esultat de P. Petersen (\cite{petersen1999},
 lemme 5.2), il existe une fonction $R(\ep)$ telle que 
$$ M= A^{R(\ep)}_{n+1},$$
avec $R(\ep) \geq \eta(\ep)$ (o\`u $\eta(\ep)$ a \'et\'e introduit 
dans la proposition \ref{a1p3}). Pour d\'emontrer la propriété (\ref{a5e19}), il suffit donc de d\'emontrer la 
propri\'et\'e de proximit\'e m\'etrique (\ref{a5e17}) pour toute 
partie $A^{R(\ep)}_k$ avec $R(\ep) \geq \eta(\ep)$.

\smallskip

\noindent {\bf Plan de la preuve de l'estimation (\ref{a5e17})}
\smallskip

\noindent Dans une premi\`ere partie, nous d\'emontrons que tout point 
de $A^{R(\ep)}_k$
(avec $R(\ep) \geq \eta(\ep)$) admet un presque antipode dans 
$A^{R(\ep)}_k$. Nous d\'emontrons \'egalement 
la proposition \ref{a9p1}. Ensuite gr\^ace au lemme \ref{a1l6}, nous montrons que 
pour tout $x,y$ dans  $A^{R(\ep)}_k$, $\cos d(x,y)$ est proche de $\sum_{i=1}^k 
a_i(x)f_i(y)$, avec $(a_i(x))_{i \geq 0}$, les coefficients de Fourier de $\cos d_x$. Nous
 montrons ensuite que $ a_i(x)$ est proche de $f_i(x)$ (il y a \'egalit\'e dans le cas de la sph\`ere), ce qui permet de 
conclure.

\subsubsection{Propri\'et\'es des ensembles $A^{R(\ep)}_k$}

Soit $k$ dans $\{2,\ldots,n+1\}$ et $A^{R(\ep)}_k$ avec
\begin{equation}\label{a5e18}
 R(\ep)\geq \eta(\ep),
\end{equation}
fix\'e.

 Sur la sph{\`e}re canonique, l'ensemble $ \{ x \in \mathbb{S}^n~;
X_1^2(x)+\dots+X_k^2(x)=1 \}$ o{\`u}  les fonctions $ (X_i)_{1\leq 
i\leq k} $ sont les $k$
premi{\`e}res fonctions coordonn\'ees, est un \'equateur de
dimension k-1. La fonction $ X_1^2+\dots+X_k^2 $
d{\'e}finie sur $ \mathbb{S}^n$ atteint son maximum sur cet 
\'equateur, son
gradient est donc nul sur cet ensemble. \\
Le lemme suivant est une g{\'e}n{\'e}ralisation de ce fait au 
{ \og presque
{\'e}quateur \fg } $ A_k^{R(\ep)}$, dans le cas o{\`u} la vari{\'e}t{\'e} 
$(M,g)$ admet $k$ valeurs propres de $n$.

\begin{lem}\label{a1l13}Il existe une fonction  $ \tau
(\epsilon) $, telle  que 
pour toute fonction $\theta(\ep)$ et pour tout élément  $(M,g)$ dans 
$\ens $ v{\'e}rifiant $ \lambda_k(M)
  \leq n+ \epsilon $, on a pour tout $x$ dans $A_k^{\theta (\epsilon)}$
 l'estimation  
$$ | \nabla (\sum_{i=1}^k f_{i}^2) | (x) \leq
4(1+\theta(\ep)) (\tau(\epsilon)+ \theta(\ep)), $$
avec $(f_i)_{1\leq i\leq k}$ une famille orthogonale de fonctions 
propres associ\'ees \`a \\
$(\lambda_i(M))_{1 \leq i \leq k}$ et 
normalis\'ees par (\ref{a9e2}).
\end{lem}
\begin{dem}
Fixons un point $x_0$ de $ A_k^{\theta(\ep)}$ et consid{\'e}rons les 
coefficients $
a_i=\frac{f_i(x_0)}{\sqrt{\sum_{i=1}^k f^2_i(x_0)}} $ pour $i$ dans 
$\{1,\ldots,k\}$. On note 
$ f=\sum_{i=1}^k a_if_i$. En d{\'e}veloppant le  terme $|\nabla f|^2 $, on obtient  
$$   |\nabla f|^2(x)= \sum_{1 \leq i,j \leq k } a_ia_j <\nabla
f_i(x),\nabla f_j (x)>. $$

\noindent D'o\`u, en $x_0$ 
$$ |\nabla f|^2(x_0)=\frac{1}{4\sum_{i=1}^k f^2_i(x_0)}|\nabla
(\sum_{i=1}^k f^2_i)(x_0)|^2.$$

\noindent Par cons\'equent en appliquant la proposition \ref{a1l3} \`a $f$ au 
point $x_0$, on obtient 
$$ \sum_{i=1}^k f^2_i(x_0) + \frac{1}{4\sum_{i=1}^k f^2_i(x_0)}|\nabla
(\sum_{i=1}^k f^2_i)(x_0)|^2 \leq 1+ \tau(\ep). $$
Or $ x_0$ appartient à $A_k^{\theta(\ep)} $ entraine $ 1+\theta(\ep)\geq 
\sum_{i=1}^k f^2_i(x_0) \geq
1-\theta(\ep)$, d'o{\`u} le r{\'e}sultat. 
\end{dem}
\smallskip

Ce lemme permet de d{\'e}montrer que tout point de $ A_k^{R(\ep)} $ 
admet un presque antipode.
\begin{lem}\label{a1l14} Il existe une fonction $ \delta(\epsilon) $ 
(v{\'e}rifiant $
  \lim_{\ep \rightarrow 0} \frac{\delta(\ep)}{\ep}=+ \infty $) telle que pour tout élément  
 $ (M,g)$ de $\ens $ v{\'e}rifiant $ \lambda_k(M)\leq n+ 
\epsilon $ et pour tout $x$ dans $ A_k^{R (\epsilon)}$, il existe
 $y$ dans $A_k^{R(\ep)}$ avec
$$    d(x,y)> \pi - \delta(\epsilon). $$
\end{lem}

\begin{dem}
Dans le cas de la sph\`ere, le point antipodal d'un point $X$ de 
$\mathbb{S}^n$ est $-X$, ce qui sugg\`ere le {\og candidat \fg } \`a \^etre 
un presque antipode de $x$ appartenant à $A^{R(\ep)}_k$. Soit $ x$ dans $A_k^{R (\epsilon)}$, notons $ \alpha = 
||(f_1,\ldots,f_k)(x)||_{\mathbb{R}^k} $, l'hypoth{\`e}se $ x$ appartient à $A_k^{R 
(\epsilon)}$
implique 
$$ (1- R(\epsilon))^{\frac{1}{2}} \leq \alpha \leq
 (1+ R(\epsilon))^{\frac{1}{2}}. $$
\noindent Par la proposition  \ref{a1p3} de presque surjectivit{\'e}, 
il existe $y$ appartenant à $A^{\eta(\ep)}_k$ tel que
$$   ||(f_1,\ldots,f_k)(y) + \frac{1}{\alpha}(f_1,\ldots,f_k)(x)|| \leq
\tau (\epsilon). $$
En appliquant le lemme \ref{a1l1} aux fonctions $(f_i)_{1 
\leq i\leq k}$, on en d\'eduit l'existence de  $ \tilde{x}$
et $\tilde{y}$ v\'erifiant 
\begin{equation}\label{a9e6} 
d(x,\tilde{x})\leq r(\epsilon),
d(y,\tilde{y})\leq r(\epsilon)
\end{equation}

\noindent et en notant $\gamma_{\tilde{x}\tilde{y}}$ l'unique 
g\'eod\'esique minimisante reliant $\tilde{x}$ \`a $\tilde{y}$, on a 
pour tout $i$ dans $ \{1,\ldots ,k\}$,
\begin{equation}\label{a1e6.6} 
\int_0^{d(\tilde{x},\tilde{y})} |(f_i \circ \gamma)''(t)+(f_i \circ 
\gamma) 
(t)|^2 dt \leq \tau'(\ep).
\end{equation}

\noindent D'autre part, par le lemme \ref{a1l13}
\begin{equation}\label{a1e6.7}
|\nabla (f_1^2+\dots+f_k^2)|(\tilde{x}) \leq \tau''(\epsilon). 
\end{equation}
Notons $ a_i=f_i(\tilde{x})$, $
b_i=(f_i\circ\gamma_{\tilde{x}\tilde{y}})'(0)$ et 
$l=d(\tilde{x},\tilde{y})$. On d\'eduit de l'{\'e}quation (\ref{a1e6.6}) et du lemme 
\ref{a1l9}, l'existence d'une fonction $\tau_{2}(\ep)$
telle que pour tout $i$ dans $\{1,\ldots,k\}$ et pour tout $t$ dans $[0,l]$,
\begin{equation}\label{a1e6.8}
 |(f_i \circ \gamma_{\tilde{x}\tilde{y}})(t)-(a_i \cos t + b_i \sin 
t) |
\leq \tau_2(\epsilon), 
\end{equation}
$$|(f_i \circ\gamma_{\tilde{x}\tilde{y}})'(t)-(-a_i\sin t + b_i\cos
t)| \leq \tau_2(\epsilon). $$
En appliquant l'inégalité de Cauchy-Schwartz, on obtient {\`a} l'aide de
(\ref{a1e6.7})
\begin{equation}\label{a8e4}
|\sum_{i=1}^{k} a_ib_i|= |\frac{1}{2}((\sum_{i=1}^k f_i^2) \circ 
\gamma_{\tilde{x}\tilde{y}})'(0)|\leq   \frac{\tau'' (\epsilon)}{2}. 
\end{equation}

Estimons maintenant $f_i(\tilde{y})$ pour $i$ dans $\{1,\ldots,k\}$.

\smallskip

\noindent 
$f_i(\tilde{y})+a_i=(f_i(\tilde{y})-f_i(y))+(f_i(y)+\frac{f_i(x)}{\alpha})+(-\frac{f_i(x)}{\alpha}+f_i(x))+(-f_i(x)+f_i(\tilde{x})).$
\noindent D'o{\`u}, par d\'efinition de $y$ et par (\ref{a9e6})
\begin{multline*} |f_i(\tilde{y})+a_i| \leq 
2C(n)r(\ep)+\tau (\ep) \\+C(n)\max
\{-(1+R(\ep))^{- \frac{1}{2}}+1;-1+(1-R(\ep))^{- \frac{1}{2}}\}.
\end{multline*}
Donc, il existe une fonction $\tau_3(\ep)$ telle que pour tout $i$ dans $\{1,\ldots,k\}$,
\begin{equation}\label{a9e7}
 |f_i(\tilde{y}) + a_i | \leq 
\tau_3(\ep).
\end{equation}
En appliquant (\ref{a1e6.8}) avec $t=l$, on obtient

\begin{equation}\label{a5e20}
 a_i \cos l + b_i \sin l = -a_i + \delta_i, 
\end{equation}
avec $ |\delta_i| \leq  \tau_3(\epsilon) $.

\noindent En multipliant (\ref{a5e20}) par $b_i$ et en sommant par rapport 
{\`a} $i$, on obtient
$$ (\sum_{i=1}^k a_ib_i)\cos l + (\sum_{i=1}^k b^2_i)\sin l  +
\sum_{i=1}^k a_ib_i \leq  (\sum_{i=1}^k b^2_i)^{\frac{1}{2}}\delta' 
.$$
avec $\delta'= \sqrt{\sum_{i=1}^k \delta_i^2}$. Or $ |b_i| \leq |\nabla f_i|$, donc la proposition 
\ref{a1l3} implique que $\sum_{i=1}^k b^2_i$ est born{\'e}e par une 
constante 
$C(n)$. D'autre part $ |\delta'|  \leq 
(n+1)^{\frac{1}{2}}\tau_3(\epsilon)$, donc on d\'eduit de 
(\ref{a8e4}) l'existence d'une 
fonction $\tau_4(\epsilon)$ telle que
$$ \left|(\sum_{i=1}^k b^2_i)\sin l\right| \leq \tau_4(\epsilon), $$
par conséquent, soit $|\sin l| \leq 
(\tau_4(\epsilon))^{\frac{1}{2}}$, soit $ \sum_{i=1}^k b^2_i\leq 
(\tau_4(\epsilon))^{\frac{1}{2}}$. 

\noindent-Premier cas : $|\sin l| \leq (\tau_4(\epsilon))^{\frac{1}{2}}$.

\noindent  Comme $||(f_1,\ldots,f_k)(x)||_{\mathbb{R}^k}$ est proche de $1$, on en 
d\'eduit que $(f_1,\ldots,f_k)(x)$ est presque égal à $-(f_1,\ldots,f_k)(y)$. Le gradient des fonctions propres étant born\'e (\ref{a9e4}), cela implique l'existence d'une constante $C'(n)>0$ telle que 
$$ d(x,y) > C'(n)$$
et donc $|\sin l| \leq (\tau_4(\epsilon))^{\frac{1}{2}}$  implique que $l$ est presque égal à $\pi$.

\noindent-Deuxi\`eme cas : $ \sum_{i=1}^k b^2_i\leq 
(\tau_4(\epsilon))^{\frac{1}{2}}
 $.

\noindent Dans ce cas par (\ref{a1e6.8}), il  existe une fonction $ 
\tau_5(\epsilon)$ telle que pour tout $i$ dans $\{1,\ldots,k\}$ et pour tout $t$ dans $[0,l]$,
$$ |(f_i \circ 
\gamma_{\tilde{x}\tilde{y}})(t))-a_i \cos t| \leq
 \tau_5(\epsilon), $$
  en appliquant cette formule avec $ t=l$, on
 obtient par (\ref{a9e7}) que  $\cos l$ est presque égal à  $-1$,  
d'o{\`u} le r{\'e}sultat.

\end{dem}

\medskip
\`A l'aide d'une l\'eg\`ere modification de la preuve ci-dessus, nous 
sommes en mesure de d\'emontrer 
la propri\'et\'e de {\og presque convexit\'e \fg} de l'ensemble 
$A^{\eta(\ep)}_k$ (proposition \ref{a9p1}).

\begin{dem} Le d\'ebut de la preuve est identique \`a celle du lemme 
\ref{a1l14}. Soit $x,y$ dans $A^{\eta(\ep)}_k$. En appliquant le lemme \ref{a1l1} aux fonctions $(f_i)_{1 \leq i \leq k}$, on en d\'eduit 
l'existence de  $ \tilde{x}$
et $\tilde{y}$ v\'erifiant $d(x,\tilde{x})\leq r(\epsilon)$, 
$d(y,\tilde{y})\leq r(\epsilon)$
 tels que, si on note $\gamma_{\tilde{x}\tilde{y}}$, l'unique 
g\'eod\'esique minimisante reliant $\tilde{x}$ \`a $\tilde{y}$, on a 
pour tout $i$ dans $\{1,\ldots ,k\}$,
\begin{equation} 
\int_0^{d(\tilde{x},\tilde{y})} |(f_i \circ \gamma)''(t)+(f_i \circ 
\gamma) (t)|^2 dt \leq 
\tau'(\ep).
\end{equation}
Notons $ a_i=f_i(\tilde{x})$, $
b_i=(f_i\circ\gamma_{\tilde{x}\tilde{y}})'(0)$ et 
$l=d(\tilde{x},\tilde{y})$. On en d\'eduit comme pr\'ecedemment (\ref{a8e4}) que pour tout $t$ dans $[0,l]$,
\begin{equation}\label{a5e25}
 |(f_i \circ \gamma_{\tilde{x}\tilde{y}})(t))-(a_i \cos t + b_i \sin 
t) |
\leq \tau_2(\epsilon) 
\end{equation}
avec
\begin{equation}\label{a5e26}
 |\sum_{i=1}^{k} a_ib_i| \leq   \frac{\tau (\epsilon)}{2}. 
\end{equation}
\noindent Par (\ref{a9e4}), il existe une constante $C(n)$ telle que 
$$ ||\nabla (\sum_{i=1}^k f_i^2) ||_{L^{\infty}} \leq C(n).$$
Par conséquent, les hypothèses sur $x$ et $\tilde{x}$ impliquent 
$$ d(x,\tilde{x}) = d_{A^{\eta(\ep)+C(n)r(\ep)}_k}(x,\tilde{x}),$$
de m\^eme pour $y$ et $\tilde{y}$. Notons 
$\eta_2(\ep)=\eta(\ep)+C(n)r(\ep)$. Pour d\'emontrer la proposition,
 il suffit donc de prouver qu'il existe des fonctions
 $\eta'(\ep)$ et $\psi(\ep)$ telles que 
$$ d_{A^{\eta'(\ep)}_k}(\tilde{x},\tilde{y})\leq 
d(\tilde{x},\tilde{y})+\psi(\ep).$$
Par construction $\tilde{y}$ appartient à $A^{\eta_2(\ep)}_k$, c'est \`a dire 
\begin{equation}\label{a5e27} 
\left| \sum_{i=1}^k (f_i^2\circ \gamma_{\tilde{x}\tilde{y}})(l)-1
\right| \leq  \eta_2(\ep).
\end{equation}
Par ailleurs, gr\^ace \`a (\ref{a5e25}) et (\ref{a5e26}), il
existe une fonction $\tau(\ep)$ telle que 
\begin{equation}\label{a5e28}
\left|\sum_{i=1}^k (f_i^2\circ \gamma_{\tilde{x}\tilde{y}})(l)-\left(
    \left(\sum_{i=1}^k a_i^2 \right) \cos^2 l+ \left(\sum_{i=1}^k
    b_i^2 \right)\sin^2 l\right)\right| 
\leq \tau(\ep).
\end{equation}
Or, comme $a_i=f_i(\tilde{x})$ et $\tilde{x}$ appartient à $A^{\eta_2(\ep)}_k$,
(\ref{a5e27}) et (\ref{a5e28}) impliquent l'existence d'une fonction
$\tau_2(\ep)$ telle que 
$$ \left|\left(1-\sum_{i=1}^k b_i^2 \right)\sin^2 l \right| \leq
  \tau_2(\ep).$$
Par cons\'equent, soit $|\sum_{i=1}^k b_i^2-1| \leq \sqrt{\tau_2(\ep)}$, soit 
$\sin^2 l \leq \sqrt{\tau_2(\ep)}$. Supposons tout d'abord que $\left|\sum_{i=1}^k b_i^2-1 \right|  \leq 
\sqrt{ \tau_2(\ep)}$. Dans ce cas, on obtient en utilisant (\ref{a5e25}) et 
(\ref{a5e26}), l'existence d'une fonction $\tau_3(\ep)$ telle que 
pour tout $t$ dans $[0,l]$, 
$$\left|\sum_{i=1}^k (f_i^2\circ 
\gamma_{\tilde{x}\tilde{y}})(t) -1 \right| \leq  \tau_3(\ep).$$
La proposition est d\'emontr\'ee dans ce  premier cas. Supposons maintenant que $\sin^2 l \leq \sqrt{ \tau_2(\ep)}$. Ce qui signifie que
$d(\tilde{x},\tilde{y})$ est proche de $0$ ou de $\pi$. Si $d(\tilde{x},\tilde{y})$ est proche de 
 $0$, c'est immédiat. 
Si  $d(\tilde{x},\tilde{y})$ est presque égal à $\pi $ alors 
n\'ecessairement, avec les notations de la proposition \ref{a1p2}, il existe $i_0$ dans 
$\{1,\ldots,k\}$  tel que $ \sin^2 d(x_{i_0},\tilde{x}) \geq \sqrt{ \tau_2(\ep)}$ et $\sin^2 d(x_{i_0},\tilde{y}) \geq \sqrt{ \tau_2(\ep)}.$ Par cons\'equent, d'apr\`es le premier cas, la courbe $c$ form\'ee de
l'union des deux g\'eod\'esiques minimisantes reliant $\tilde{x}$ \`a
 un point voisin de $x_{i_0}$ et ce dernier point \`a $\tilde{y}$ est 
contenue dans $
A^{\tau_3(\ep)}_k$. Le lemme \ref{a1l2} sur la fonction 
{ \og excess \fg} permet alors de conclure.
\end{dem}

\subsubsection{D\'emonstration de la propri\'et\'e de {\og proximit\'e 
m\'etrique \fg}}

Nous venons de montrer que tout point de $A^{R(\ep)}_k$ admet un 
presque antipode. L'id\'ee de la preuve consiste \`a utiliser la 
propri\'et\'e des fonctions $\cos d_p$ pour $p$ admettant un presque 
antipode, \'etablie dans le lemme \ref{a1l6}. Cependant, on ne peut 
pas appliquer directement le lemme \ref{a1l6} avec la fonction 
$\delta(\ep)$ introduite dans
le lemme \ref{a1l14}, puisqu'on voit facilement (par exemple dans la
proposition \ref{a1l4}) que les fonctions $\tau (\ep)$ utilis\'ees 
dans la
preuve du lemme \ref{a1l14} sont
sup{\'e}rieures {\`a} $\ep$ et donc 
\begin{equation}\label{rajout1}
\lim_{\ep \rightarrow 0} \frac{\delta (\ep)}{\ep}= + \infty
\end{equation} 
 Pour contourner ce problème, on pose
$$ \ov[k]= \max \{ i~; \lambda_i \leq n+ \sqrt{\delta (\ep)} \}.$$
Par (\ref{rajout1}) et pour $\ep$ assez petit, on a  $ \ov[k] \geq k$.
 D'apr{\`e}s le lemme \ref{a1l6} appliqu\'e avec 
$\eta=\delta(\ep)$ et $\sqrt{\delta(\ep)}$ et le lemme \ref{a1l14}, 
il existe une 
fonction $\tau(\ep)$ telle que pour tout $ x$ dans $A_k^{R (\epsilon)}$, 
il existe des coefficients $(\alpha_i
(x))_{i=1}^{\ov[k]} $ pour lesquels pour tout $z$ dans $M$, 
\begin{equation}\label{a1e6.9}
 | \cos d_x(z)- \sum_{i=1}^{\ov[k]} \alpha_i (x)f_i(z) | \leq 
\tau
(\epsilon), 
\end{equation}
avec $ | \sum_{i=1}^{\ov[k]} {\alpha}^2_i (x)-1| \leq \tau (\epsilon).$

Montrons que $  \sum_{i=1}^{\ov[k]} | f_i(x)-{\alpha}_i(x)|^2 $ est 
petit. 
$$ \sum_{i=1}^{\ov[k]} | f_i(x)-{\alpha}_i(x)|^2 =  
\sum_{i=1}^{\ov[k]}  f^2_i (x) + 
 \sum_{i=1}^{\ov[k]} {\alpha}^2_i (x) -2  \sum_{i=1}^{\ov[k]} 
{\alpha}_i (x) f_i(x)
 .$$
Or le lemme \ref{a1l11} implique 
\begin{equation}\label{a5e29}
  \sum_{i=1}^{\ov[k]}  f^2_i (x)
\leq 1+ \tau (\sqrt{\delta (\epsilon)})
\end{equation}
D'autre part, l'inégalité (\ref{a1e6.9}) appliquée pour $z=x$ donne
$$ | \sum_{i=1}^{\ov[k]} {\alpha}_i (x) f_i(x) -1| \leq \tau 
(\epsilon) .$$
On obtient finalement
\begin{equation}\label{a1e6.91}
  \sum_{i=1}^{\ov[k]} | f_i(x)-{\alpha}_i(x)|^2 \leq 3 \tau 
(\epsilon) +
 \tau (\sqrt{\delta (\epsilon)}). 
\end{equation}
\noindent En appliquant l'inégalité (\ref{a1e6.9}) à $y$ dans $A^{R(\ep)}_k$, il 
vient 
$$ \left| < F(x),F (y)>_{\mathbb{R}^k}- \cos d(x,y) \right| \leq 
\left| \sum_{i=1}^{k} f_i(x)f_i(y) -
  \sum_{i=1}^{\ov[k]} {\alpha}_i(x)f_i(y) \right| + \tau(\ep)$$
$$ \leq \left| \sum_{i=1}^{\ov[k]}
  \left(f_i(x)-{\alpha}_i(x)\right)f_i(y) \right| + \left|
  \sum_{i=k+1}^{\ov[k]}f_i(x)f_i(y) \right| + \tau (\ep). $$
En utilisant l'in\'egalit\'e de Cauchy-Schwartz puis (\ref{a1e6.91}) 
et (\ref{a5e29}), on 
obtient 
$$ \left| \sum_{i=1}^{\ov[k]}
  \left(f_i(x)-{\alpha}_i(x)\right)f_i(y) \right| \leq \left( \left( 
3\tau
  (\ep) + \tau ( \sqrt{\delta (\ep)}) \right) \left( 1+
  \tau (\sqrt{\delta(\ep)}) \right)
\right)^{\frac{1}{2}}. $$
Il ne reste plus qu'{\`a} estimer le terme $ |
\sum_{i=k+1}^{\ov[k]}f_i(x)f_i(y)| $. Or pour tout {\'e}l{\'e}ment z de $A_k^{R(\ep)}$, on a
$$ \sum_{i=1}^k f_i^2 (z) \geq 1- R(\ep),$$
donc en utilisant de nouveau l'in\'egalit\'e de Cauchy-Schwartz et 
(\ref{a5e29}), on obtient \\
$$ |
\sum_{i=k+1}^{\ov[k]}f_i(x)f_i(y)| \leq R(\ep) +  \tau
  (\sqrt{\delta(\ep)}), $$
ce qui termine la démonstration.


\end{document}